\documentclass[11pt]{amsart}
\usepackage{graphicx}
\newtheorem{thm}{Theorem}[section]

\newtheorem{prop}[thm]{Proposition}
\theoremstyle{definition}
\newtheorem{defn}[thm]{Definition}
\theoremstyle{remark}
\newtheorem{rem}[thm]{Remark}
\numberwithin{equation}{section}

\newcommand{\CC}{\mathbb C}
\newcommand{\ZZ}{\mathbb{Z}}

\newcommand{\C}{\mathcal{C}}

\begin{document}

\title{Real Surfaces in Elliptic Surfaces}%
\author{Marko Slapar}%
\address{Institute of Mathematics, Physics and Mechanics, University of Ljubljana, Jadranska 19, 1000 Ljubljana, Slovenia}%
\email{marko.slapar@fmf.uni-lj.si}%
\date{June 27, 2003}

\begin{abstract}
We study the structure of complex points on real surfaces,
embedded into complex Elliptic surfaces. We show, for example,
that any compact surface has a totally real embedding into a
blow-up of a $K3$ surface. We also exhibit smooth disc bundles
over compact orientable surfaces that have a Stein structure as
Stein domains inside Elliptic surfaces.
\end{abstract}

\maketitle

\section{Statement of results}

\par Let $S$ be a real surface, embedded into a complex surface $(X,J)$.
We say that the embedding $S\hookrightarrow X$ is \emph{totally
real} at a point $p\in S$, if $T_pS+J(T_pS)=T_pX$. If this is not
the case, we call $p$ a \emph{complex point} of the embedding. An
embedding is called totally real, if it is totally real at all
points.
\begin{thm} \label{1}Every compact oriented real surface $S$ has a totally real
embedding into any $K3$ surface. Every compact real surface has a
totally real embedding into a blow-up of a $K3$ surface at one
point.
\end{thm}
A blow-up of a complex surface is of course not minimal. If we
want to have an embedding of all compact surfaces into a minimal
surface, we have the following theorem.
\begin{thm}\label{2} Every compact real surface has a totally real embedding into
any $E(3)$ surface.
\end{thm}
Let us denote by $\Sigma_g$ the compact Riemann surface of genus
$g\ge 0$, and let $n$ be an integer. We denote by $D(g,n)$ the
open unit disc bundle over the surface $\Sigma_g$, with Euler
number $n$. It follows from the adjunction inequality for Stein
surfaces that for $n>2g-2$, the smooth manifolds $D(g,n)$ do not
have any Stein structure. It is furthermore a consequence of the
result of Gompf \cite{Gompf1}, using a method of Stein surgery
developed by Eliashberg in \cite{Eli1}, that for $n\le 2g-2$, the
smooth manifolds $D(g,n)$ can be endowed with a Stein structure.
We use a method of Stein fattening, introduced by Forstneri\v c
\cite{For1}, to give a different proof of this result, by
explicitly seeing $D(n,g)$ as open strictly pseudoconvex Stein
domains in Elliptic surfaces $E(n)$. The definition of Elliptic
surfaces $E(n)$ is given latter in the text.
\begin{thm} \label{3}For $n\le 2g-2$, the smooth manifold $D(n,g)$ can be
realized as an open strictly pseudoxconvex Stein domain in a
complex Elliptic surface $E(-n+2g)$.
\end{thm}
For compact nonorientable surfaces, the situation is even less
rigid as in Theorem \ref{3}. Let us denote $\tilde D(n,\chi)$ the
disc bundle over a compact nonorientable surface with Euler
characteristic $\chi$ and Euler number $n$. By passing to a double
cover and using adjunction inequality for Stein surfaces, $\tilde
D(n,\chi)$ can have a Stein structure only if $n+\chi\le 0$. See
\cite{Gompf0}. We can construct such Stein structures using the
following theorem.
\begin{thm} \label{4} For $n+\chi\le 0$, the smooth manifold $\tilde D(n,\chi)$ can be
realized as an open strictly pseudoxconvex Stein domain in
sufficiently high blow-up of $\CC P^2$.
\end{thm}
\section{Complex points of real surfaces in complex surfaces}
Let $S$ be a real compact surface in a complex surface $X$. After
perhaps a small generic perturbation, the embedding
$S\hookrightarrow X$ is totally real outside a finite collection
of complex points, which can, following Bishop \cite{Bis}, be
classified as either \emph{elliptic} or \emph{hyperbolic}. This
means that, locally around  a complex point $p$, we can choose
complex coordinates $(z,w)$ on $X$, so that $S$ is written as
$$w=\alpha z\bar {z} + {1\over 2}z^2 + {1\over 2}\bar{z}^2+o(|z|^3),$$
where $\alpha\in [0,\infty]$ is a holomorphic invariant of the
complex point (the case $\alpha=\infty$ should be understood as
the surface $w=z\bar z+o(|z|^3)$). Elliptic points correspond to
$\alpha>1$, hyperbolic to $\alpha<1$ and (nongeneric) parabolic to
$\alpha=1$.
\par Let $S\hookrightarrow X$ be a generic embedding of a compact surface with only
finitely many complex points, all of them either elliptic or
hyperbolic. Let $e(S)$ be the number of elliptic complex points on
$S$, and $h(S)$ a number of hyperbolic complex points on $S$. The
algebraic number of complex points, $I(S)=e(S)-h(S)$, is a
topological invariant of the embedding, and can be expressed as
\begin{equation}\label{Web}I(S)=\chi(S)+\chi(NS),\end{equation}
where $\chi(S)$ is the Euler characteristic of $S$ and $\chi(NS)$
is the Euler characteristic of the normal bundle $NS$ of the
embedding, which can be calculated as the self-intersection number
$[S]^2$ in the case of an orientable $S$. For the proof, see
\cite{Web}. In the case of an oriented surface $S\hookrightarrow
X$, one can also talk about a sign at a complex point $p\in S$:
positive, if the orientation of the tangent space $T_pS$ agrees
with the induced orientation on $T_pS$ as a complex subspace of
$T_pX$, and negative otherwise. The signed algebraic numbers of
complex points, $I_\pm(S)=e_\pm (S)-h_\pm (S)$, also turn out to
be topological invariants. Here the sign in the subscripts
indicate the sign of complex points. They can be expressed by Lai
formulae \cite{Lai}
\begin{equation}\label{Lai} 2I_\pm(S)=\chi(S)\pm
\left<c_1(X),[S]\right>+[S]^2,\end{equation} where
$c_1(X)=c_1(TX)$ is the first Chern class of the complex surface
$X$.
\par By a result of Harlamov and Eliashberg \cite{ElH}, and more
generally Forstneric \cite {For1}, a pair of elliptic and
hyperbolic complex points can always be cancelled on nonoriatable
surfaces, and can be cancelled on orientable surfaces exactly when
the signs at the complex points agree. The cancellation can be
done by a $\C^0$ small isotopy, changing the embedding only in an
arbitrary small neighborhood of an arc, connecting the two complex
points.
\par By a result of Bishop \cite{Bis}, a surface $S\hookrightarrow X$ having elliptic
complex points can never have a Stein neighborhood basis in $X$,
since the local hull of holomorphy is nontrivial at elliptic
complex points. Up to changing an embedding by an isotopy, this is
the only obstruction to having a Stein basis: if $S$ is a compact
oriented surface in a complex surface $X$ with $I_\pm(S)\le 0$,
then $S$ is $\C^0$ isotopic to a surface having a regular Stein
neighborhood basis, see \cite{For1}. The same result holds if the
surfaces $S$ is compact and nonoriantable, with the topological
condition now being $I(S)\le 0$. This is most easily proven by
first cancelling all the elliptic points by an isotopy, and than
by further isotoping the surface, putting the remaining hyperbolic
points to be of a special type, so that the usual distance
function gives a regular Stein neighborhood basis. For a careful
proof, see \cite{For1,For2}. By a regular basis of $S$, we mean a
basis system of open neighborhoods
$\{U_\epsilon\}_{0\le\epsilon\le 1}$ of $S$ in $X$, satisfying
\begin{itemize}
\item $\Omega_\epsilon=\bigcup_{s< \epsilon}\Omega_s$, \item
$\bar\Omega_\epsilon=\bigcap_{s>\epsilon}\Omega_s$, \item
$S=\bigcap_{s>0}\Omega_s$ is a strong deformation retract of
$\Omega_\epsilon$. \end{itemize}
\par Up to first performing a small isotopy of a surface, finding compact surfaces
with regular Stein neighborhood basis reduces to checking whether
$I_\pm$ (or $I$ if the surface is unorientable) is nonpositive for
the given embedding. It turns out that this is very often
automatically the case. Using Seiberg-Witten theory,
Kronheimer-Mrowka \cite{KrM}, Fintushel-Stern \cite{FiS} and
Oszvath-Szabo \cite{OzS}, have proved adjunction inequalities for
surfaces in many $4$-manifolds, implying that $I_\pm(S)\le 0$ for
any oriented surface imbedded in a compact K\"{a}hler surface with
$b_2^+>1$, except for the embedded spheres, where one has a
slightly weaker inequality, $I_\pm\le 2$. This led to the result
of Lisca-Mati\'{c} \cite{LiM}, stating that $I_\pm (S)\le 0$ holds
for all oriented real surfaces embedded into Stein surfaces,
except for homologically trivial spheres. We call this result the
adjunction inequality for Stein surfaces.
\section{Elliptic surfaces }
\begin{defn} A complex elliptic surface is a compact complex surface $X$,
together with a holomorphic map $f\! :X\to C,$ where $C$ is a
compact complex curve, so that all but finitely many fibers
$f^{-1}(z)$ are elliptic curves.
\end{defn}
Here we only study a special kind of Elliptic surfaces, called
$E(n)$ surfaces. Up to diffeomorphism, they can be given by the
following construction.
\par Let $\CC P^2 \rightarrow \CC P^1$ be a pencil of cubics over $\CC
P^1$. By this, we mean a collection of curves
$\{z_0p_0+z_1p_1,[z_0,z_1]\in \CC P^1\}$ in $\CC P^2$, where
$p_0,p_1$ are generic (intersecting at $9$ distinct points) cubics
in $\CC P^2$. Let $E(1)=\CC P^2 \# 9\overline{\CC P^2}\rightarrow
\CC P^1$ be the fibration gotten by blowing up at the $9$
singularities. The generic fiber of this fibration is an elliptic
curve, but since $\chi(E(1))=12$, we must also have some singular
fibres. To construct surfaces $E(n)$ we need the fibre sum
operation. Let us assume we have constructed $E(1),...,E(n-1)$.
Let $F_1$ be a regular fibre of $E(1)\to \CC P^1$ and let $F_2$ be
a regular fibre of $E(n-1)\to \CC P^1$. Let $U_1$ and $U_2$ be
tubular neighborhoods of $F_1$ and $F_2$ respectively, and let
$\phi\!:\partial U_1 \to
\partial U_2$ be fibre preserving, orientation reversing
diffeomorphism. Then $E(n)=E(1)\#_fE(n-1):=(E(1)\backslash
U_1\bigcup E(n-1)\backslash U_2)/\phi$.
\par We gave a definition of $E(n)$, usually found in the
literature. It is very convenient for calculating properties of
$E(n)$, but from it, it is not obvious that the diffeomorphism in
the fibre sum operation can be chosen so the surfaces $E(n)$ have
a complex structure. This is indeed the case, as one can see from
alternative descriptions of surfaces $E(n)$. For a more algebraic
treatment of Elliptic surfaces, see for example \cite{GH}.
\par We now list some homological properties of
Elliptic surfaces $E(n)$. The results are classical, and can be
found in \cite{Gompf0}.
\begin{prop} The $4$-manifold $E(n)$
is simply connected with $H_2(E(n)),\ZZ)=\ZZ^{12n-2}$ and
intersection form
$$n(-E_8)\oplus 2(n-1)\left[ \begin{array}{rr}
      0 & 1 \\ 1 & -2
             \end{array} \right]\oplus
    \left[ \begin{array}{rr}
      0 & 1 \\ 1 & -n
             \end{array} \right],
$$
where all basis elements of $H_2(E(n),\ZZ)$ are represented by
embedded spheres, except the elements with self intersection $0$,
which can be represented by embedded by tori. The homology element
with self intersection $-n$ is the class of a section (sphere) of
the fibration.
\end{prop}\begin{prop} For any $E(n)$ surface we have $c_1(E_n)=(2-n)PD(f)$,
where $PD(f)$ is the Poincare dual of the class of a fiber of the
elliptic fibration $E(n)\to\CC P^1$.
\end{prop}
Since $c_1(E(2))=0$, and $E(2)$ is simply connected, $E(2)$
surfaces are exactly $K3$ surfaces.
\begin{rem} It turns out that a minimal simply connected Elliptic
surfaces with sections is always diffeomorphic to an $E(n)$
surface with $n>1$, see \cite{Kas}. \end{rem}
\section{Proofs of theorems}
\begin{proof}[Proof of theorem \ref{1}.] Let $X$ be any $K3$
surface. In the homology group $H_2(X,\ZZ)$, take $s$ and $f$ to
be basis homology classes, generating one of the $3$ summands with
intersection matrix $\left[\begin{array}{rr}
      0 & 1 \\ 1 & -2
             \end{array}\right]$.
Let $S$ be a sphere with homology class $s$ and $F_1,\cdots F_g$
be $g$ nonintersecting tori in the homology class $f$, all
intersecting the sphere $S$ at exactly one positive transverse
intersection. Let $\Sigma$ be the surface, gotten by resolving all
intersections of the union $S\cup F_1\cup\cdots\cup F_g$. By a
resolution of intersections, we simply mean substituting a local
model $zw=0$ of an intersection by $zw=\epsilon$, where $\epsilon$
is small. Then $\Sigma$ is an oriented surface of genus $g$ in the
homology class $s+gf$. Applying formulas $(\ref{Lai})$, we get
$2I_\pm(\Sigma)=\chi(\Sigma)+[\Sigma]^2=0$. So after a small
isotopic perturbation, we can make $\Sigma$ totally real. This way
we can construct $3$ non-intersecting, non-homologous totally real
embeddings of a compact oriented surface of genus $g$ into $X$.
\par Let us now look at embeddings of nonorientable surfaces. Let $\Sigma$ be a compact
nonorientable surface with Euler characteristics
$\chi(\Sigma)=\chi$. By a result of Massey \cite{Mas}, $\Sigma$
can be embedded into $\CC^2$ with $\chi(N\Sigma)\in
\{2\chi-4,2\chi,\ldots,4-2\chi\}$. By the same result, this are
the only possible normal Euler numbers for such an embedding.
Using $(\ref{Web})$, we can achieve $I(\Sigma)$ to be any one of
the numbers in the set
$$N(\chi):=\{3\chi-4,3\chi,\ldots,4-\chi\},$$
by embedding $\Sigma$ into $\CC^2$. \\
\textbf{Case $1$:} $\chi\equiv 0 (\textrm{mod}\ 4)$. Since $0$ is
in the set $N(\chi)$, we can embed $\Sigma$ in a small
contractible domain inside any complex surface.\\
\textbf{Case $2$:} $\chi\equiv 2 (\textrm{mod}\ 4)$. Since $2$ is
in $N(\chi)$, we can embed $\Sigma$ into a small contractible set
in $K3$ with $I(\Sigma)=2$. Let $S$ be a totally real sphere in
$K3$, not intersecting $\Sigma$, and let $\Sigma'=\Sigma\# S$. The
connected sum is performed by simply tubing an arc between the
surfaces. Then $\Sigma'$ is again a nonorientable surface with
Euler characteristic $\chi$. Calculating the algebraic number of
complex points using $(\ref{Web})$, we have
$I(\Sigma')=I(\Sigma)+I(S)-2=0$. \\
\textbf{Case $3$:} $\chi\equiv 3 (\textrm{mod}\ 4)$. Let $E$ be
the exceptional sphere in the blow-up of $K3$. Then
$I(E)=\chi(E)+[E]^2=1.$ Since $1\in N(\chi)$, let us embed
$\Sigma$ into a small contractible set in $K3$, disjoint with $E$,
so that $I(\Sigma)=1$. Let $\Sigma'=\Sigma\# E$. As before,
$\Sigma'$ is nonorientable with Euler characteristic $\chi$. We
have $I(\Sigma')=I(\Sigma)+I(E)-2=0$.\\
\textbf{Case $4$:}$\chi\equiv 1 (\textrm{mod}\ 4)$. Let $E$ be the
exceptional sphere in the blow-up of $K3$ and $S$ be the totally
real sphere in $K3$. We assume they are disjoint. Since $3\in
N(\chi)$, let us embed $\Sigma$ into a small contractible set in
$K3$, disjoint with $E$ and $S$, so that $I(\Sigma)=3$. Let
$\Sigma'=\Sigma\# E\# S$. $\Sigma'$ is nonorientable with Euler
characteristic $\chi$ and $I(\Sigma')=I(\Sigma)+I(E)+I(S)-4=0$.
This completes the proof.
\end{proof}
\begin{proof}[Proof of theorem \ref{2}.] Since, as above, $H_2(E(3),\ZZ)$ group contains a
pair of a sphere and a torus, spanning a direct summand in
homology with intersection matrix $\left[\begin{array}{rr}
      0 & 1 \\ 1 & -2
             \end{array}\right]$,
we only need to review Cases $3$ and $4$ in the above proof. So
let $\Sigma$ be a nonorinetable surface with
$\chi(\Sigma)=\chi\equiv 3(\textrm{mod}\ 4)$. We can embed
$\Sigma$ into $E(3)$ having $I(\Sigma)=5$. Let $S$ be a totally
real sphere in $E(3)$ and let $S'$ be a disjoint section of
$E(3)$. We have $I(S')=\chi(S')+[S']^2=-1$. The surface
$\Sigma'=\Sigma\# S\# S'$ is nonorientable with
$\chi(\Sigma')=\chi$ and $I(\Sigma')=I(\Sigma)+I(S)+I(S')-4=0$. A
small perturbation produces a totally real embedding. If
$\chi(\Sigma)=\chi\equiv 1(\textrm{mod}\ 4)$, we first embed
$\Sigma$ into $E(3)$ with $I(\Sigma)=3$ and then make
$\Sigma'=\Sigma\# S'$.
\end{proof}
\begin{proof}[Proof of theorem \ref{3}] Let $S$ be a section of
the fibrtion $E(m)\rightarrow \CC P^1$ and let $F_1,\ldots,F_g$ be
$g$ generic fibers. Let $\Sigma=S\bigcup_{i=1}^gF_i$ and let
$\tilde \Sigma$ be a smooth surface we get from resolving the
intersections of $\Sigma$. Of course, $[\Sigma]=[\tilde \Sigma]$
and $\chi(\tilde \Sigma)=2-2g$. By using $c_1(E(m))=(2-m)PD(f)$,
where $f$ is the homology of a fiber, and applying  the formula
\ref{Lai}, we get $I(\tilde\Sigma)=0$ and $I_+(\tilde\Sigma)=2-m$.
We can thus perturb $\tilde\Sigma$ to a surface, having regular
Stein neighborhood basis.  Since
$[\tilde\Sigma]^2=([S]+gf)^2=-m+2g$, we  must set $2g-n=m$, so
that the elements of the Stein neighborhood basis are
diffeomorphic to $D(n,g)$.
\end{proof}
\begin{proof}[Proof of theorem \ref{4}] Let $\Sigma$ be a nonorientable compact surface with
$\chi(\Sigma)=\chi$. We use a similar argument as in the proof of
Theorem \ref{1}. As there, we can embed $\Sigma$ into a small
contractible set in $\CC P^2$ so that $\chi(N\Sigma)\in
\{2\chi-4,2\chi,\ldots,4-2\chi\}$. Let $\CC P^2\#m\overline {\CC
P^2}$ be the blow-up of $\CC P^2$ at $m$ distinct points, and let
$E_1,\ldots,E_m$ be the exceptional spheres. Then
$\Sigma'=\Sigma\# E_1\cdots\#E_m$ is again a nonorientable surface
with Euler characteristic $\chi$ and
$\chi(N\Sigma')=\chi(N\Sigma)-m$. For a right choice of $m$, we
can always achieve that $n\in \{2\chi-4,2\chi,\ldots,4-2\chi\}-m$.
We can thus choose an embedding of $\Sigma$ so that $\Sigma'$ has
$\chi(N\Sigma')=n$. Since $I(\Sigma')=\chi+n\le 0$, a perturbation
of $\Sigma'$ has a regular strictly pseudoconvex Stein
neighborhood basis of the right topological type.
\end{proof}
\begin{rem} Instead of a blow-up of $\CC P^2$, we could also use
surfaces $E(m)$ to put Stein structures on $\tilde D(n,\chi)$.
Instead of taking connected sums with $m$ exceptional spheres in
the blow-up, we take just one connected sum with a section of
$E(m)$.
\end{rem}
\subsection*{Acknowledgment} I would like to thank Franc
Forstneri\v c, for many useful discussions on the topic.


\begin{thebibliography}{}
\bibitem{Bis}
E. Bishop, Differential Manifolds in complex Euclidian space, Duke
Math. J. \textbf{32} (1965), 1-22

\bibitem{Eli1}
Y. Eliashberg, Topological characterization of Stein manifolds of
dimension $>2$, Internat. J. of Math. \textbf{1} (1990), 29-46
\bibitem{ElH}
 Y. Eliashberg and V.M. Harlamov, On the number of complex points of a real surface in a
complex surface, Proc. Leningrad int. Topology Conf.
(1982),143-148
\bibitem{FiS}
R. Fintushel and R. Stern, Immersed spheres in $4$-manifolds and
the immersed Thom conjecture, Turk. J. Math. \textbf{19} (1995),
145-157
\bibitem{For1}
F. Forstneri{\v c}, Complex tangents of real surfaces in complex
spaces, Duke Math. J. \textbf{67} (1992), 353-376
\bibitem{For2}
 F. Forstneri{\v c}, Stein domains in complex surfaces, J. Geom. Anal. \textbf{13} (2003), 77-94
\bibitem{Gompf0}
R. Gompf and A. Stipsicz, $4-$Manifolds and Kirby Calculus,
Graduate Studies in Math. \textbf{20}, A.M.S., 1999
\bibitem{Gompf1}
R. Gompf, Handlebody construction of Stein surfaces, Ann. of
Math.(2) \textbf{148} (1998), 619-693
\bibitem{GH}
P. Griffits and J. Harris, Principles of algebraic geometry,
Willey, 1978
\bibitem{Kas}
A. Kas, On the deformation types of regular elliptic surfaces,
Complex Analysis and Algebraic geometry (1977), 107-111
\bibitem{KrM}
P.B. Kronheimer and T. Mrowka, The genus of embedded surfaces on
the projective plane, Math. Res. Lett. \textbf{1} (1994), 797-808
\bibitem{Lai}
H. F. Lai, Characteristic classes of real manifolds immersed in
complex manifolds, Trans. AMS \textbf{172} (1972), 1-33
\bibitem{LiM}
P. Lisca and G. Matic, Tight contact structures and Seiberg-Witten
invariants, Invent. Math. \textbf{128} (1997), 509-525
\bibitem{Mas}
W. S. Massey, Proof of a conjecture of Whitney, Pacific J. Math.
\textbf{31} (1969), 143-156
\bibitem{Nem}
S. Nemirovski, Complex analysis and differential topology on
complex surfaces, Uspekhi Math. Nauk \textbf{45} (1999), 47-74
\bibitem{OzS}
P. Ozsv{\'{a}}th and Z. Szab{\'{o}}, The symplectic Thom
conjecture, Ann. of Math. (2) \textbf{151} (2000), 93-124
\bibitem{Web}
S.M. Webster, Minimal surfaces in a K{\"{a}}hler surface, J. Diff.
Geom. \textbf{20} (1984), 463-470
\end{thebibliography}
\end{document}